# Most Probable Dynamics of Stochastic Dynamical Systems with Exponentially Light Jump Fluctuations

Yang Li[1, 2, a], Jinqiao Duan[2, b], Xianbin Liu[1, c] and Yanxia Zhang[3, 2, d]

[1] State Key Laboratory of Mechanics and Control of Mechanical Structures, College of Aerospace Engineering, Nanjing University of Aeronautics and Astronautics, 29 Yudao Street, Nanjing 210016, China

[2] Department of Applied Mathematics, Illinois Institute of Technology, Chicago, Illinois 60616, USA

[3] Department of Mechanics, Beijing Institute of Technology, Beijing 100081, China

[a] li_yang@nuaa.edu.cn

[b] duan@iit.edu

[c] Corresponding author: xbliu@nuaa.edu.cn

[d] zhangyanxia1314@126.com

**Abstract:** The emergence of the exit events from a bounded domain containing a stable fixed point induced by non-Gaussian Lévy fluctuations plays a pivotal role in practical physical systems. In the limit of weak noise, we develop a Hamiltonian formalism under the Lévy fluctuations with exponentially light jumps for one- and two-dimensional stochastic dynamical systems. This formalism is based on a recently proved large deviation principle for dynamical systems under non-Gaussian Lévy perturbations. We demonstrate how to compute the most probable exit path and the quasi-potential by several examples. Meanwhile, we explore the impacts of the jump measure on the quasi-potential quantitatively and on the most probable exit path qualitatively. Results show that the quasi-potential can be well estimated by an approximate analytical expression. Moreover, we discover that although the most probable exit paths are analogous to the Gaussian case for the isotropic noise, the anisotropic noise leads to significant changes of the structure of the exit paths. These findings shed light on the underlying qualitative mechanism and quantitative feature of the exit phenomenon induced by non-Gaussian noise.

**Keywords:** Rare events; transition under non-Gaussian fluctuations; Hamiltonian formalism; most probable exit path; action functional and quasi-potential

**The exit phenomenon from a domain of attraction for a stable state in weak noise limit plays a fundamental role in a broad range of nonlinear systems. Based on Hamiltonian formalism or path integral formulation, this problem has been well studied for systems with (Gaussian) Brownian noise. However, a large number of dynamical systems with bursting or intermittent features arise in various scientific fields. It is more appropriate to model these systems with non-Gaussian Lévy noise. To uncover the internal mechanisms of these rare events in exit phenomena driven by non-Gaussian noise is of great importance. We establish a Hamiltonian formalism for stochastic differential equations with exponentially light jump fluctuations, based on a large deviation principle. Meanwhile, the numerical experiments reveal the impacts of the jump measure on the most probable exit path, and also on quasi-potential which provides the information about how and when the exit events occur.**

## 1. Introduction

Environmental noisy fluctuations are inevitable in dynamical systems and may lead to unexpected or rare events. Even weak background noise may have a profound impact in dynamics if the observations are performed on a sufficiently long time scale. For instance, it may result in the transition phenomenon between distinct metastable states which can be observed in various scientific fields, such as biology, physics, chemistry, engineering and finance[1-6]. To explore the mechanism of transition behavior is a challenging task in stochastic dynamical systems.

The last several decades have witnessed a great deal of mathematical and experimental efforts devoted to the research of rare events with the methods of Hamiltonian formalism or equivalent path integral formulations[7-10]. Under weak random perturbations, Freidlin and Wentzell[11] established the mathematical foundation of large deviation principle. The associated action functional describes the possibility of the stochastic path passing through the vicinity of a given trajectory. The theory asserts that the rare events, when they occur, will follow a certain trajectory with an overwhelming probability. This specific trajectory is referred to as the most probable exit path, or the optimal exit path. Along this trajectory the action functional attains its global minimum which is defined as the quasi-potential. Its value exponentially dominates the magnitude of the expected exit time. Therefore, the two crucial quantities, the most probable exit path and the quasi-potential, are theoretically investigated based on a variational principle.

Large deviation theory for dynamical systems under Gaussian noise has been extensively applied to the dynamical systems with various structures such as fixed points[12, 13], limit cycles[14] and chaotic attractors[15, 16]. Smelyanskiy et al.[14] investigated the patterns of extreme paths in the vicinity of an unstable focus of a periodically oscillating system. Maier and Stein[17] studied the Gaussian noise-driven exit events in a double-well system lacking detailed balance and found a bifurcation of the most probable exit path as a parameter varies. Chen et al.[16] revealed that the escape from a nonhyperbolic chaotic attractor occurs through a hierarchical sequence of crossings between stable and unstable manifolds of saddle cycles.

However, certain complex phenomena are not suitable to be modeled as stochastic differential equations with Gaussian noise, due to peculiar dynamical features such as the abrupt events in climate system[18-20] and the burst-like events in the transcription process of gene regulation[4, 21, 22]. A stochastic process with discontinuous trajectories, i.e., non-Gaussian Lévy process appears more appropriate in describing the fluctuations in these systems. The exit problem for non-Gaussian Lévy noise is still under development because of the complicated nonlocal term in the statistical distribution of the noisy fluctuations.[23, 24] Several authors investigated the exit phenomena in neurosystems under non-Gaussian Lévy noise to disclose its excitation behavior.[25, 26, 27] Gao et al.[28] designed an efficient and accurate numerical scheme to compute the mean exit time and escape probability for stochastic differential equations with non-Gaussian Lévy noise. Zheng et al.[19] developed a probabilistic framework to investigate the maximum likelihood climate change for an energy balance system under the influence of greenhouse effect and non-Gaussian $\alpha$-stable Lévy motions. Gomes[29] proved that the solutions of stochastic differential equations with exponentially light jump Lévy fluctuations satisfy a large deviation principle and then estimated the magnitude of its exit time. Imkeller et al.[30] made further efforts to estimate the order of the exponential pre-factor of the mean first exit time for this type process. Based on their works, we attempt to deduce the Hamiltonian formalism and locate the initial conditions for the exit paths of dynamical systems under the exponentially light jump fluctuations. Then we compute the most probable exit path and the quasi-potential, and compare the similarities and differences with the Gaussian case.

This article is organized as follows. We describe the Hamiltonian formalism and provide an example of the exit problem for one- and two-dimensional dynamical systems in Sections 2 and 3, respectively. Then conclusions are drawn in Section 4.

## 2. One-dimensional stochastic dynamical systems

### 2.1. Hamiltonian formalism

Consider a one-dimensional stochastic dynamical system

$$dx_t = f(x_t)dt + dL_t^\varepsilon, \tag{1}$$

where $L_t^\varepsilon = \varepsilon L_{t/\varepsilon}$ with $L_t$ being a Lévy process with exponentially light jumps (see Appendix). The jump measure $\nu$ satisfies $\nu(dy) = e^{-\lambda|y|^\alpha} dy$ for $\alpha > 1$ and $\lambda > 0$ [29, 30]. We like to know when the small noise intensity $\varepsilon$ approaches zero, the most probable exit path and expected exit time from a certain bounded domain containing a stable state of the underlying deterministic dynamical system.

According to Refs.[29, 31-33], the solutions of Eq. (1) satisfy a large deviation principle with action functional

$$S_{t_0 t}[\varphi(s)] = \int_{t_0}^t L(\varphi(s), \dot{\varphi}(s))ds \tag{2}$$

in the sense that

$$P_{x_0}\left\{\sup_{t_0 \le s \le t}|x(s) - \varphi(s)| < \delta\right\} \sim \exp\left\{-\varepsilon^{-1} S_{t_0 t}[\varphi(s)]\right\}, \tag{3}$$

where $\varphi(t)$ is an absolutely continuous function and the subscript $x_0$ denotes the common initial point $x(t_0) = x_0$ and $\varphi(t_0) = x_0$. Here $L$ may be regarded as a Lagrangian function. It implies that the probability of the sample trajectories passing through a small tube of a certain function is exponentially estimated by the action functional (2). Thus the action functional actually dominates the possibility that a trajectory realizes. Then its global minimum corresponds to the path with largest probability. Therefrom, the Euler-Lagrange equation, whose solution provides the most probable exit path, can be obtained by variation of the action functional as

$$\frac{d}{dt}\left(\frac{\partial L}{\partial \dot{\varphi}}\right) - \frac{\partial L}{\partial \varphi} = 0. \tag{4}$$

However, the Lagrangian in Eq. (2) cannot be expressed explicitly as in the Gaussian case. Fortunately, the corresponding Hamiltonian, which is the Legendre transform of the Lagrangian, does have the following analytical expression[32]

$$H(x, p) = f(x)p + \int_{R\setminus 0}\left(e^{yp} - 1 - yp\chi_{|y|\le 1}\right)\nu(dy). \tag{5}$$

Here, $p = \partial L/\partial \dot{x}$ is referred as the momentum and $R$ denotes the set of real numbers. In fact, the integration of the term $yp\chi_{|y|\leq 1}$ vanishes since the measure $\nu$ is symmetric about $y$. According to classical results of analytical mechanics, the Euler-Lagrange equation can be transformed into the equivalent auxiliary Hamiltonian system

$$\dot{x} = \frac{\partial H}{\partial p} = f(x) + \int_{R\setminus 0} y e^{yp-\lambda|y|^\alpha} dy,$$
$$\dot{p} = -\frac{\partial H}{\partial x} = -f'(x)p.$$
(6)

Therefore, the complicated computation for the probability is transformed into the integral solution of this Hamiltonian system and this is the reason why the method is called as Hamiltonian formalism. The projection $x(t)$ of the solution $(x(t), p(t))$ to the coordinate space provides the most probable exit path.

In order to quantify the expected exit time, Freidlin and Wentzell[11] defined an effective and novel concept, i.e., the quasi-potential, also called nonequilibrium potential

$$W(x) = \inf_{t_0, t} \inf_{\varphi(t_0)=x_1, \varphi(t)=x} S_{t_0 t}[\varphi(s)],$$
(7)

where $x_1$ is the stable fixed point of the underlying deterministic system. In other words, the value of the quasi-potential at $x$ is defined as the minimum of the action functional throughout all the possible paths connecting $x$ and the fixed point. As we can infer from its name, it measures the relative stability of separate metastable state in the sense that the expected exit time from the domain of attraction of the stable state exponentially depends on the depth of the potential well of quasi-potential

$$E\tau_{exit} \sim \exp\left(\varepsilon^{-1} \inf_{x\in\partial U_{x_1}} W(x)\right),$$
(8)

where $\partial U_{x_1}$ denotes the boundary of basin of the stable state $x_1$. Therefore, the exit phenomenon can be essentially identified as the transition of distinct potential wells through overcoming the barrier from the perspective of energy. Note that $H(x,0)=0$ from Eq. (5). Hence, $H(x_1,0)=0$. Due to the facts that the most probable path connects the fixed points and that Hamiltonian remains constant along the continuous connecting trajectories, $H(x,p)=0$ throughout the most probable path. Remark that the most probable path beginning at the fixed point $x_1$ requires $t_0 \to -\infty$ as $x \to x_1$. Therefore, $W(x(t)) = \int_{-\infty}^{t} \dot{x}(s)\cdot p(s) ds$ along the most probable path, or equivalently,

$$\dot{W} = \dot{x}\cdot p = f(x)p + \int_{R\setminus 0} yp e^{yp-\lambda|y|^\alpha} dy.$$
(9)

Note that Eqs. (6) and (9) construct a complete group of differential-integral equations to evolve the exit path and quasi-potential. However, we still need initial conditions to perform numerical integration. Theoretically, the most probable path should start at the fixed point as the initial time approaches negative infinity while it is impossible in practical terms. Consequently, the initial points should be chosen at the unstable manifold of the fixed point and close to the fixed point in the extended phase space, as shown in Fig. 2(a). Therefore, we need to carry out local linearization around the fixed point. Denote $B$ as the Jacobian matrix of Eq. (6)

$$B = \begin{bmatrix} C & D \\ 0 & -C \end{bmatrix}, C \equiv f'(x_1), D \equiv \int_{R \backslash 0} y^2 e^{-\lambda |y|^\alpha} dy \tag{10}$$

and it is seen that the matrix $B$ has two eigenvalues $C<0$ and $-C>0$. Then the unstable eigenvector of the fixed point corresponding to the eigenvalue $-C$ is $(1, -2C/D)^T$. Consequently, the initial condition for the most probable exit path is selected as

$$x_0 = x_1 + \delta, \ p_0 = -2C\delta/D, W_0 = -C\delta^2/D, \tag{11}$$

where $\delta$ is a small parameter.

**2.2 Example: A stochastic energy balance model**

Unlike the Gaussian case, theoretical computation is difficult for the exit behavior even in the one-dimensional system due to the nonlocal integral term in the Hamiltonian formalism. Therefore, we take a specific model into consideration to perform numerical calculations on the basis of the discussions in Section 2.1.

The system chosen for investigation is a stochastic energy balance model describing the climate change[19]

$$dT = -U'(T)dt + dL_t^\varepsilon \tag{12}$$

with the potential function

$$U(T) = \frac{1}{Ch}\left(-\frac{1}{4}S_0\left(0.5T + 2\ln\left(\cosh\frac{T-265}{10}\right)\right) + \frac{1}{5}\gamma\theta T^5\right), \tag{13}$$

where $T$ represents the global mean surface temperature and the heat capacity $Ch = 46.8 \text{Wyrm}^{-2}$ is defined as the amount of heat that must be added to the object in order to raise its temperature. The solar constant $S_0 = 1368 \text{Wm}^{-2}$, the Stefan constant $\theta = 5.67 \times 10^{-8} \text{Wm}^{-2}\text{K}^{-4}$ and the greenhouse factor $\gamma \in [0,1]$. The parameter $\gamma = 0.61$ corresponds to the bistable case with two stable fixed points $x_1$, $x_2$ separated by one unstable state $x_0$, for which its potential function is shown in Fig. 1.

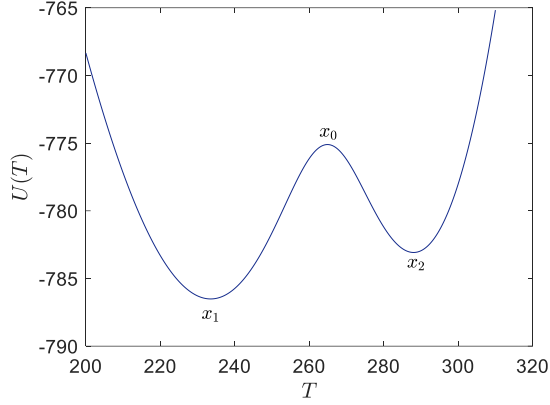

Fig. 1. The potential function of the energy balance model.

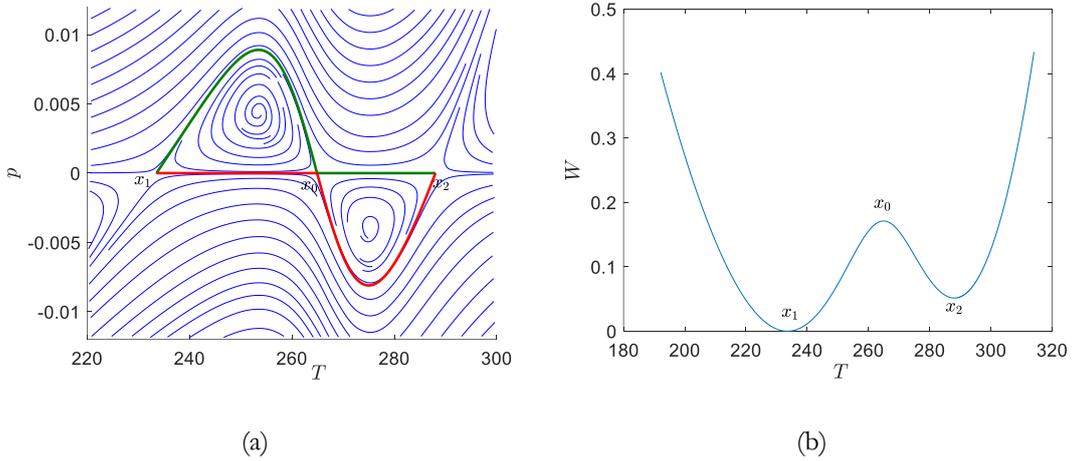

(a)　　　　　　　　　　　　　　　(b)

Fig. 2. (a) The extended Hamiltonian phase space. The green and red curves indicate the most probable transition paths from $x_1$ to $x_2$ and from $x_2$ to $x_1$, respectively. (b) Quasi-potential landscape.

Based on the specific dynamical system, we integrate Eqs. (6) and (9) combined with the initial condition (11) to evolve the exit path and quasi-potential and plot them in Fig. 2. The green curve indicates the most probable path from $x_1$ to $x_2$ and the red curve is the reverse in Fig. 2(a). Note that although the sample paths of the noise are discontinuous, the most probable path is smooth sufficiently. Reviewing the auxiliary Hamiltonian system (6), we can find that the corresponding deterministic system is recovered if the momentum is set to be zero. In other words, although it is not so apparent as in the Gaussian case, the momentum also takes the place of noise and can be regarded as an indicator of the effect of noise. It is seen that the momentum is small within the vicinity of the fixed point (stable or unstable) and is large in other area before escaping. After the particle passes over the unstable saddle state, the momentum is exactly zero. This is because the particle will follow the deterministic flow to approach the other stable state so that it does not cost extra energy during this process.

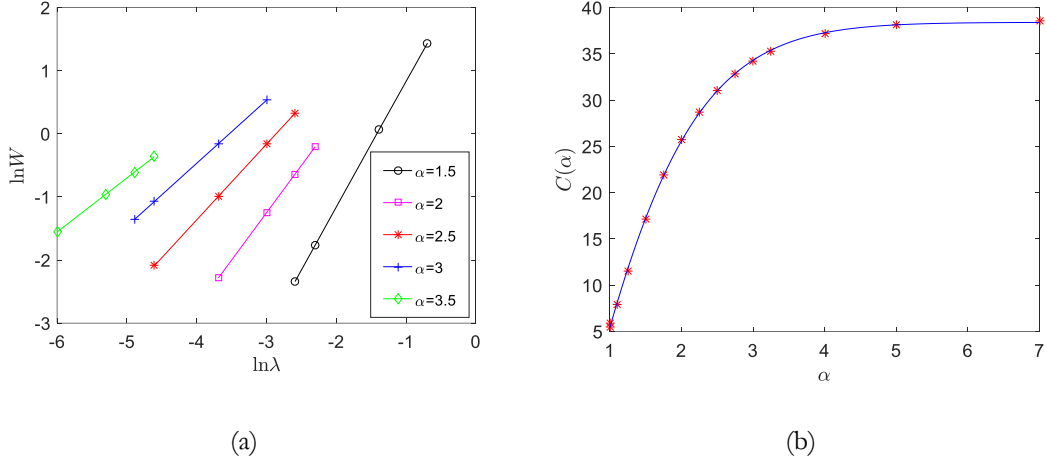

(a)                        (b)

Fig. 3. (a) The logarithm of quasi-potential versus the one of $\lambda$. (b) The function between $C(\alpha)$ and $\alpha$. The red stars indicate the numerical results and the blue curve is the fitting function.

Table 1: Relation between $\varphi(\alpha)$ and $\alpha$

| $\alpha$ | 1.5 | 2 | 2.5 | 3 | 3.5 |
|---|---|---|---|---|---|
| $\varphi(\alpha)$ | 1.9913 | 1.4991 | 1.1989 | 0.9987 | 0.8567 |

In addition, it is found that the shape of the quasi-potential in Fig. 2(b) is exactly the same as the potential function of the model except its magnitude. This phenomenon is similar to the case with Gaussian noise in which the quasi-potential is precisely twice the potential function, although the present noise is much more complicated. However, its magnitude is significantly changed and affected by the Lévy jump measure. Now we identify the depth of the quasi-potential well of $x_1$ as a representative quantity to show its dependency relations about the parameters $\lambda$ and $\alpha$.

First, we observe the relationship between the quasi-potential and $\lambda$ with $\alpha$ fixed as $\alpha = 1.5, 2, 2.5, 3, 3.5$. Fig. 3(a) shows an obvious linear relation between $\ln W$ and $\ln \lambda$ which implies $W = C(\alpha)\lambda^{\varphi(\alpha)}$. Via the least square method, we can compute the slopes of these lines and list them in Table 1. According to these datum, it is so apparent that $\varphi(\alpha) = 3/\alpha$. Dividing $W$ by $\lambda^{3/\alpha}$, we can show the relationship between $C(\alpha)$ and $\alpha$ by the red star sign in Fig. 3(b). It seems that this function has a horizontal asymptotic line so that we can fit it by the function $C(\alpha) = k_2 - k_1 \exp(-s_2 \alpha^{s_1})$. It is seen that the fitting curve with $s_1 = 1.363$, $s_2 = 0.6009$, $k_1 = 60.21$ and $k_2 = 38.41$ almost crosses every point. Above all, the relationship between the quasi-potential and the parameters in the jump measure can be written as

$$W \approx \left[ k_2 - k_1 \exp(-s_2 \alpha^{s_1}) \right] \lambda^{3/\alpha}, \tag{14}$$

where the coefficients rely on the specific system. In order to confirm it, an arbitrary combination of the parameters such as $\lambda = 0.68$ and $\alpha = 1.83$ is chosen to calculate the quasi-potential. The result for Eq. (14) is 0.11975058 which agrees perfectly with the one for numerical integral 0.11981649. Once the quasi-potential is computed, the expected exit time can be estimated exponentially according to its value on the saddle point.

## 3. Two-dimensional stochastic dynamical systems

### 3.1. Hamiltonian formalism

Consider a two-dimensional stochastic dynamical system

$$\begin{aligned} dx_1 &= f_1(x_1, x_2) dt + dL_{1t}^\varepsilon, \\ dx_2 &= f_2(x_1, x_2) dt + dL_{2t}^\varepsilon, \end{aligned} \tag{15}$$

where $L_{1t}^\varepsilon$ and $L_{2t}^\varepsilon$ are independent scalar Lévy processes with exponentially light jumps of small noise intensity as in Section 2. Their jump measures satisfy $\nu_i(dy_i) = \exp\{-\lambda |y_i|^{\alpha_i}\} dy_i$ for $\lambda > 0$ and $\alpha_i > 1$, $i = 1, 2$. Again, we consider the exit phenomenon from the domain of attraction of the stable state $\bar{\mathbf{x}}$ of the underlying deterministic dynamical system.

Similarly, the Hamiltonian for two-dimensional system is

$$H(\mathbf{x}, \mathbf{p}) = f_1 p_1 + f_2 p_2 + \int_{R \setminus 0} \int_{R \setminus 0} \left( e^{y_1 p_1 + y_2 p_2} - 1 \right) \nu_1(dy_1) \nu_2(dy_2). \tag{16}$$

The auxiliary Hamiltonian system which evolves the exit paths is then

$$\begin{aligned}
\dot{x}_1 &= f_1(x_1, x_2) + \int_{R\backslash 0}\int_{R\backslash 0} y_1 \exp\{y_1 p_1 + y_2 p_2\} \nu_1(dy_1)\nu_2(dy_2), \\
\dot{x}_2 &= f_2(x_1, x_2) + \int_{R\backslash 0}\int_{R\backslash 0} y_2 \exp\{y_1 p_1 + y_2 p_2\} \nu_1(dy_1)\nu_2(dy_2), \\
\dot{p}_1 &= -\partial f_1/\partial x_1 \cdot p_1 - \partial f_2/\partial x_1 \cdot p_2, \\
\dot{p}_2 &= -\partial f_1/\partial x_2 \cdot p_1 - \partial f_2/\partial x_2 \cdot p_2
\end{aligned} \qquad (17)$$

and along each trajectory, the quasi-potential can be integrated by

$$\dot{W} = \dot{x}_1 \cdot p_1 + \dot{x}_2 \cdot p_2. \qquad (18)$$

In fact, the original stable fixed point is transformed into a saddle point by adding two unstable eigenvectors and the deterministic coordinate space is embedded into the four-dimensional Hamiltonian space as the two-dimensional stable manifold of the original fixed point. Furthermore, the trajectories that are described by the solutions of Eq. (17), and that originate from the fixed point as $t_0 \to -\infty$, span a two-dimensional unstable Lagrangian manifold in four-dimensional extended phase space. The projections of these paths onto the coordinate space determine the patterns of the extreme paths. By saying the extreme paths, we refer to the paths corresponding to the local minimum of the action functional. Since there may exist several extreme paths arriving at the same point, the value of quasi-potential is only provided by the global minimum of their action functional. This phenomenon is an essential difference between one- and two-dimensional systems. In order to evolve these paths, the initial conditions should be chosen at the Lagrangian manifold around the fixed point. Denote $B$ as the Jacobian matrix of Eq. (17)

$$\begin{aligned}
B &= \begin{bmatrix} C & D \\ 0 & -C^T \end{bmatrix}, C \equiv \frac{\partial \mathbf{f}}{\partial \mathbf{x}}, D \equiv \begin{bmatrix} d_{11} & d_{12} \\ d_{21} & d_{22} \end{bmatrix}, \\
d_{11} &= \int_{R\backslash 0}\int_{R\backslash 0} y_1^2 \nu_1(dy_1)\nu_2(dy_2), \, d_{22} = \int_{R\backslash 0}\int_{R\backslash 0} y_2^2 \nu_1(dy_1)\nu_2(dy_2), \\
d_{12} &= d_{21} = \int_{R\backslash 0}\int_{R\backslash 0} y_1 y_2 \nu_1(dy_1)\nu_2(dy_2).
\end{aligned} \qquad (19)$$

The two unstable eigenvectors of $B$ are indicated as $\mathbf{e}_{u1}$ and $\mathbf{e}_{u2}$. Then the point at the unstable manifold around the fixed point can be represented as

$$\begin{pmatrix} \mathbf{x} \\ \mathbf{p} \end{pmatrix} = c_{u1} \begin{pmatrix} \mathbf{e}_{u1\mathbf{x}} \\ \mathbf{e}_{u1\mathbf{p}} \end{pmatrix} + c_{u2} \begin{pmatrix} \mathbf{e}_{u2\mathbf{x}} \\ \mathbf{e}_{u2\mathbf{p}} \end{pmatrix}. \qquad (20)$$

After eliminating the coefficients, we obtain

$$\mathbf{p} = \begin{pmatrix} \mathbf{e}_{u1\mathbf{p}} & \mathbf{e}_{u2\mathbf{p}} \end{pmatrix} \begin{pmatrix} \mathbf{e}_{u1\mathbf{x}} & \mathbf{e}_{u2\mathbf{x}} \end{pmatrix}^{-1} \mathbf{x}. \qquad (21)$$

Denote $M = \begin{pmatrix} \mathbf{e}_{u1\mathbf{p}} & \mathbf{e}_{u2\mathbf{p}} \end{pmatrix} \begin{pmatrix} \mathbf{e}_{u1\mathbf{x}} & \mathbf{e}_{u2\mathbf{x}} \end{pmatrix}^{-1}$ and then the initial conditions for exit paths have the following forms

$$\delta \mathbf{x}_0 = (r\cos\theta, r\sin\theta), \ \mathbf{x}_0 = \bar{\mathbf{x}} + \delta \mathbf{x}_0,$$
$$\mathbf{p}_0 = M\delta \mathbf{x}_0, \ W_0 = \frac{1}{2}\delta \mathbf{x}_0^T M \delta \mathbf{x}_0, \tag{22}$$

where $r$ is a small fixed positive parameter and $\theta \in [0, 2\pi)$. In fact, this choice of the initial conditions can be generalized to higher dimensional dynamical systems.

### 3.2 Example: Maier-Stein system

The classical Maier-Stein system is chosen as a model for investigation[17]

$$dx_1 = \left(x_1 - x_1^3 - \gamma x_1 x_2^2\right)dt + dL_{1t}^\varepsilon,$$
$$dx_2 = -\left(1 + x_1^2\right)x_2 dt + dL_{2t}^\varepsilon, \tag{23}$$

where $\gamma$ is a positive parameter. It is easily seen that there exist two stable nodes SN1(-1,0) and SN2(1,0) whose domains of attraction are separated by the stable manifold of the unstable saddle point US(0,0), as shown in Fig. 4. Note that the vector field is symmetric about $y$-axis so that we just consider the exit phenomenon from the basin of SN1.

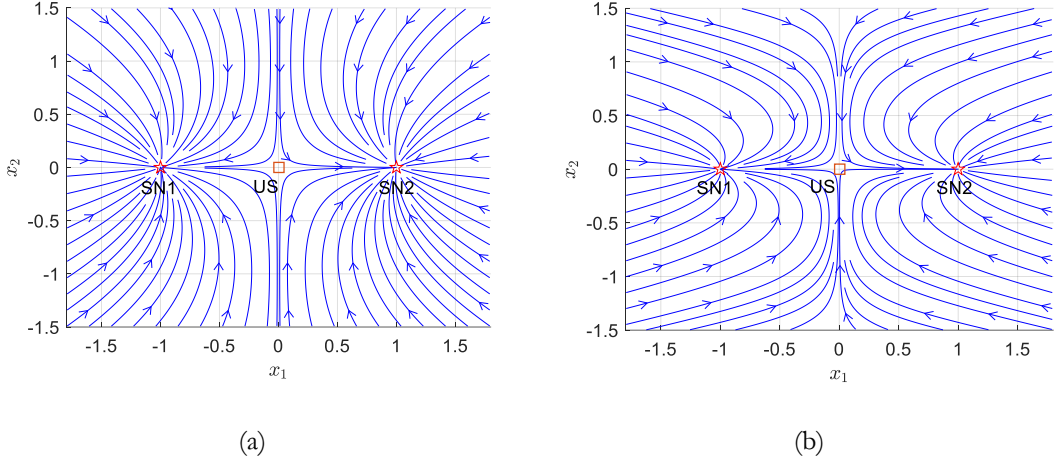

(a)          (b)

Fig. 4. (a) The vector field of Maier-Stein system for $\gamma = 1$. (b) The vector field of Maier-Stein system for $\gamma = 5$.

**Isotropic case:** $\alpha_1 = \alpha_2 = \alpha$

Under the isotropic random perturbations, the action plot method[34] can be used to provide all the actions of the extreme paths to exit the basin of the stable state. The so-called actions imply the values of action functional on the extreme paths. Their values are recorded when we integrate Eqs. (17) and (18) combined with the initial conditions (22) until the boundary. Thus these actions as well as their corresponding extreme paths are characterized by the initial angle positions. The minimal action provides the exponential magnitude of the expected exit time and its corresponding exit path appears as the most probable one. It is the quasi-potential of the saddle point.

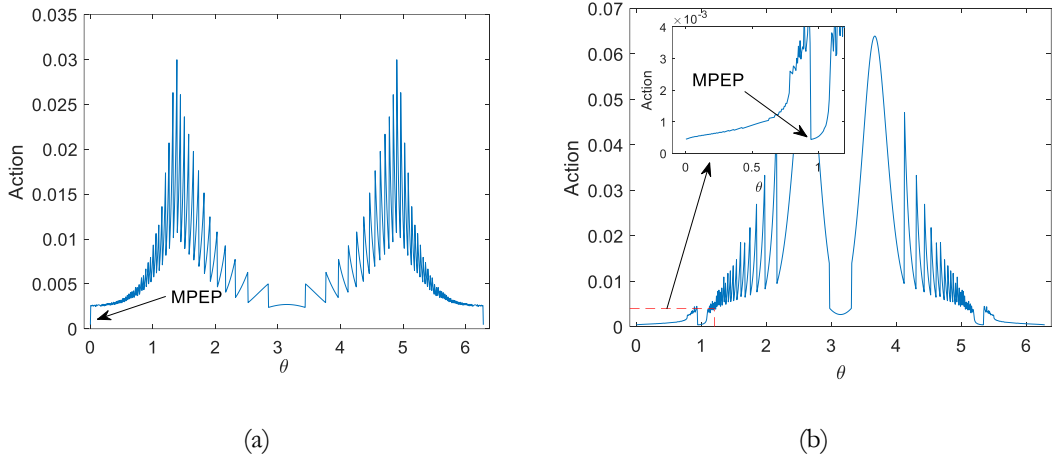

Fig. 5. (a) Action plot for $\gamma=1$. (b) Action plot for $\gamma=5$. The inset figure is local amplification of the rectangular region.

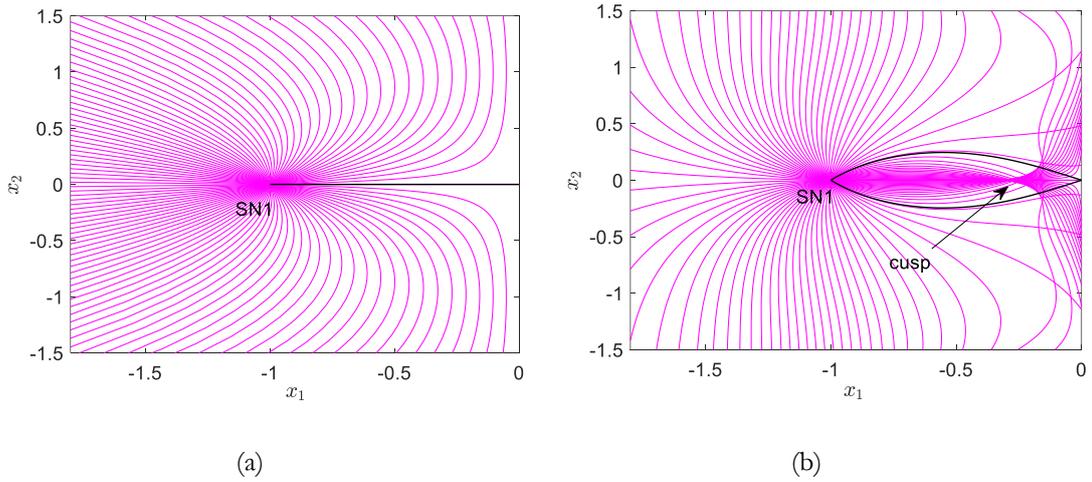

Fig. 6. (a) The patterns of extreme paths from SN1 for $\gamma=1$. The black line denotes the most probable exit path. (b) The patterns of extreme paths from SN1 for $\gamma=5$. The two black curves denote the most probable exit paths.

A large amount of initial points are selected on a small circle around the fixed point. Then the results of action plot are depicted in Fig. 5 for both $\gamma=1$ and $\gamma=5$ with the parameters $\lambda=0.1$ and $\alpha=1.5$. We depict the partial results of the extreme paths and the most probable exit paths in Fig. 6 for both $\gamma=1$ and $\gamma=5$. In the case of $\gamma=1$ corresponding to the gradient field, the most probable exit path is the horizontal line connecting SN1 and US which is just the inverse of the unstable manifold of US. However, it ceases to be the most probable one but bifurcates to two symmetric most probable paths which are the heteroclinic trajectories of SN1 and US as $\gamma = 5$. This bifurcation actually stems from the appearance of the cusp indicated in Fig. 6(b). The Lagrangian manifold in the right side of the cusp is folded into three pieces. Thus the projection of the two-dimensional Lagrangian manifold into the coordinate space is not one-to-one mapping such that there may be several extreme paths reaching one point. For example, there are three extreme paths reaching US, i.e., the two most probable exit paths and the horizontal line. However, the action of the horizontal one is larger than the others such that it is no longer the most probable. According to Refs.[17, 35], this singularity and the structures of these extreme paths for different parameters are qualitatively analogous to the Gaussian case.

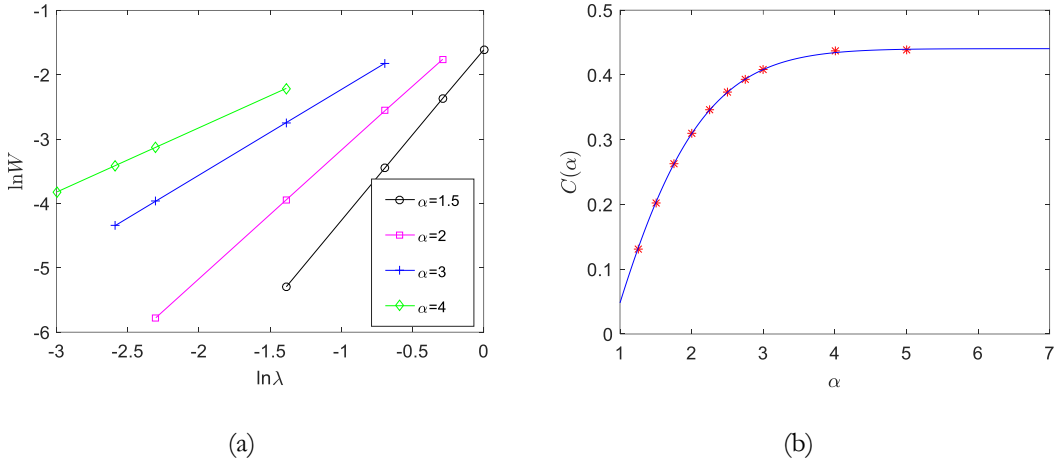

(a)                          (b)

Fig. 7. (a) The logarithm of quasi-potential versus the one of $\lambda$. (b) The function between $C(\alpha)$ and $\alpha$. The red stars indicate the numerical results and the blue curve is the fitting function.

Table 2: Relation between $\varphi(\alpha)$ and $\alpha$

| $\alpha$ | 1.5 | 2 | 3 | 4 |
|---|---|---|---|---|
| $\varphi(\alpha)$ | 2.6574 | 1.9944 | 1.3316 | 0.9991 |

However, its quasi-potential is greatly affected by the non-Gaussianity of noise. We still compute the quasi-potential of the saddle point as the one-dimensional case to show its change. Without loss of generality, the non-gradient field case of $\gamma = 5$ is taken into consideration. With $\alpha$ fixed as 1.5, 2, 3, 4 respectively, the logarithm of the quasi-potential shows obvious linear relation about the logarithm of $\lambda$ in Fig. 7(a). The slopes of these lines are evaluated and listed in Table 2. Clearly, it is found that $W = C(\alpha)\lambda^{4/\alpha}$. The function $k_2 - k_1 \exp(-s_2\alpha^{s_1})$ can be still used to fit the curve of $C(\alpha)$ versus $\alpha$. It agrees with the results of numerical integration in Fig. 7(b) when the coefficients are chosen as $s_1 = 1.487$, $s_2 = 0.6062$, $k_1 = 0.7196$ and $k_2 = 0.4405$. Above all, the quasi-potential can be represented as

$$W \approx \left[ k_2 - k_1 \exp(-s_2\alpha^{s_1}) \right] \lambda^{4/\alpha} \tag{24}$$

in two-dimensional stochastic dynamical systems.

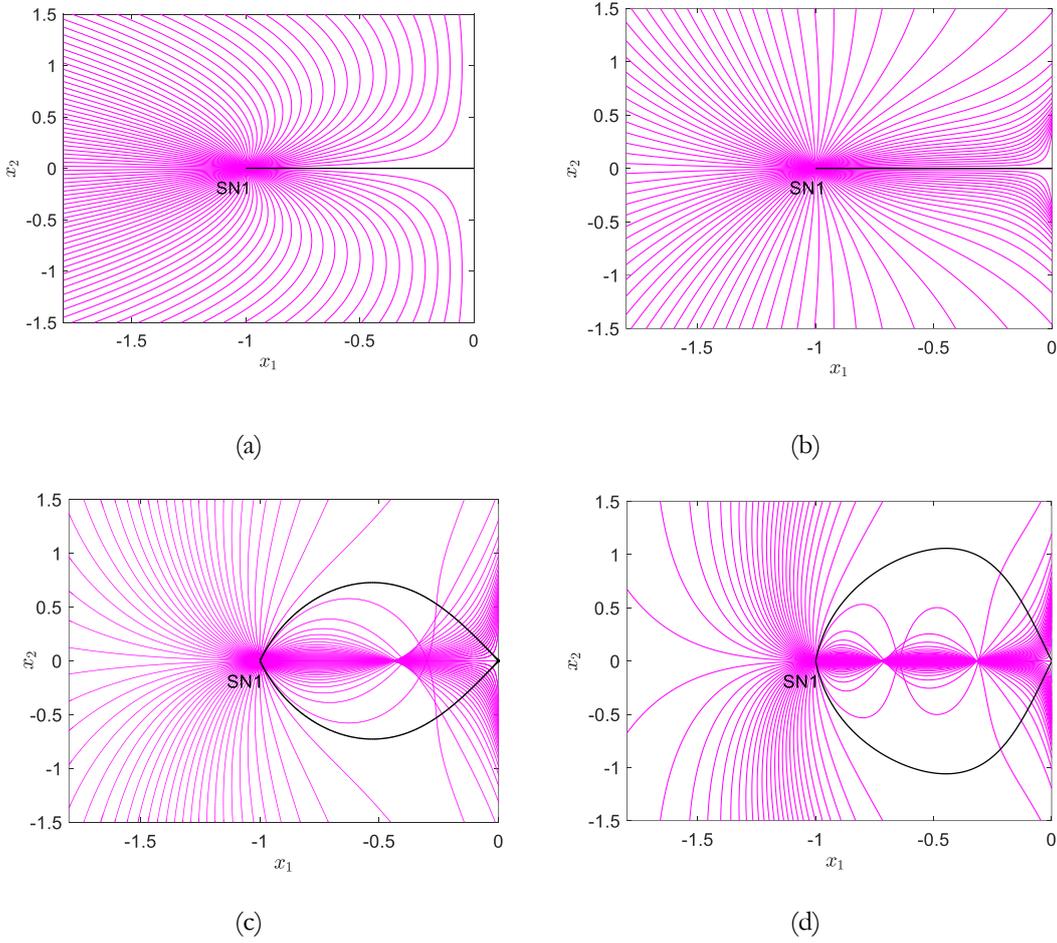

(a) (b) (c) (d)

Fig. 8. The patterns of extreme paths from SN1 for $\gamma = 1$ and $\alpha_1 = 1.5, 2, 2.5, 3.5$ respectively. The black curves denote the most probable exit paths.

**Anisotropic case:** $\alpha_1 \neq \alpha_2$

With regard to the anisotropic case, it is difficult to provide an analytical expression for the quasi-potential since too many parameters exist. Furthermore, the patterns of the extreme paths will be significantly changed if $\alpha_1 \neq \alpha_2$. A group of the parameters of $\lambda = 0.1$, $\alpha_2 = 1.5$ and several values of $\alpha_1$ is taken as an example. Fig. 8 shows the extreme paths and the most probable ones for $\alpha_1 = 1.5, 2, 2.5, 3.5$ and $\gamma = 1$. Even in this gradient field situation, the increase of $\alpha_1$ will induce the appearance of the singularity of the Lagrangian manifold. It is found that there exist one and two cusps for the cases of $\alpha_2 = 2.5$ and $\alpha_2 = 3.5$, respectively. Moreover, as denoted by the black curves, the most probable path also bifurcates from the horizontal line to two symmetric heteroclinic trajectories of SN1 and US with increasing $\alpha_1$.

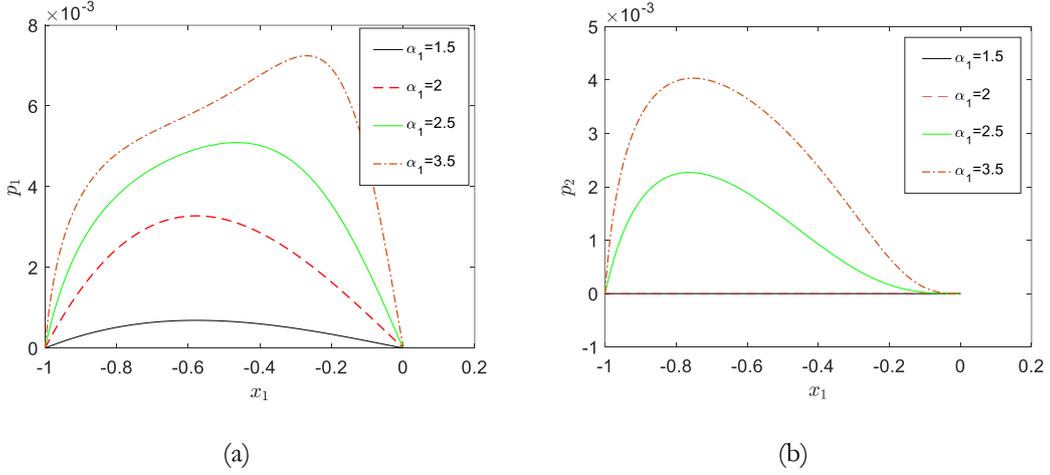

(a)  (b)

Fig. 9. (a) The optimal fluctuational force $p_1$ for various $\alpha_1$. (b) The optimal fluctuational force $p_2$ for various $\alpha_1$.

The momenta $p_1$ and $p_2$ corresponding to the optimal fluctuational forces that drive the particle to realize the most probable exit path are depicted in Fig. 9. When $\alpha_1 = 1.5$, we have $p_2 = 0$ and it is only the noise of the first direction that takes effect since the exit path is horizontal. The same behavior occurs for $\alpha_1 = 2$ while the magnitude of $p_1$ is greatly amplified. As $\alpha_1$ increases further, $p_2$ is no longer zero, which means the noise of the second direction starts to be effective and its proportion becomes larger gradually. Therefore, the increase of $\alpha_1$ actually results in weaker noise in the first direction such that it costs too much energy to exit horizontally. Consequently, the most probable exit path is bent deeper and deeper since the vertical force can play a more important role to reduce the energy needed.

## 4. Conclusion

In this paper, we have investigated the exit problem of one- and two-dimensional stochastic dynamical systems with exponentially light jump fluctuations in the weak noise limit. In both cases, we have established the Hamiltonian formalism and located the initial conditions lying in the unstable Lagrangian manifold around the fixed point to compute the most probable exit path and the quasi-potential. We presented several examples to illustrate our method.

Specifically, we have considered a stochastic energy balance model for global climate as an example in the one-dimensional case. The numerical results show an obvious analytical relationship between the quasi-potential and the parameters of the jump measure which is perfectly fitted. The shape of the quasi-potential is found the same as the potential function of the deterministic system. In addition, taking the Maier-Stein system as an example, we delineate the patterns of the extreme paths to uncover its singularity and reveal the bifurcation of the most probable exit path. In the isotropic case, the value of the quasi-potential is again estimated quantitatively. Furthermore, it is found that the extreme paths including the most probable exit path are similar to that of the Gaussian case. But in the anisotropic case, the most probable exit paths are significantly different for the Gaussian case, as the singularity of the Lagrangian manifold and the bifurcation of the most probable exit path can be induced by varying the jump measure.

In short, we have developed an efficient method to describe and compute the exit paths of stochastic dynamical systems with a specific class of non-Gaussian Lévy noise, i.e., those fluctuations with the exponentially light jumps (or abrupt bursts). This is theoretically based on a large deviation principle for dynamical systems driven by exponentially light jump fluctuations for $\alpha > 1$, in the weak noise limit[29, 32]. However, the large deviation principle is no longer valid for $0 < \alpha < 1$ but a moderate deviation principle holds[29].

Finally, we note that stochastic dynamical systems with the well-known $\alpha$-stable Lévy motions are not yet shown to satisfy a large deviation principle,[23] and further efforts are needed to examine the exit phenomenon for these non-Gaussian systems.

**Acknowledgements**

This research was supported by the National Natural Science Foundation of China (Grants 11472126 and 11232007), the Priority Academic Program Development of Jiangsu Higher Education Institutions (PAPD), and the China Scholarship Council (CSC No. 201906830018).

## Appendix. Lévy processes

A scalar Lévy process $L_t$ is a stochastic process satisfying the following conditions:

(i) $L_0 = 0$, a.s.;

(ii) Independent increments: for any choice of $n \geq 1$ and $t_0 < t_1 < \cdots < t_{n-1} < t_n$, the random variables $L_{t_0}$, $L_{t_1} - L_{t_0}$, $L_{t_2} - L_{t_1}$, $\cdots$, $L_{t_n} - L_{t_{n-1}}$ are independent;

(iii) Stationary increments: $L_t - L_s$ and $L_{t-s}$ have the same distribution;

(iv) Stochastically continuous sample paths: for every $s > 0$, $L_t \to L_s$ in probability, as $t \to s$.

The Lévy-Itô decomposition of Lévy process with Lévy triplet $(b, \sigma^2, \nu)$ is

$$L_t = bt + \sigma B_t + \int_0^t \int_{|y|>1} y N(dy, ds) + \int_0^t \int_{0<|y|\leq 1} y \tilde{N}(dy, ds),$$

where $B_t$ is a Brownian motion, $N$ denotes the jump counting measure and $\tilde{N}$ is the compensated jump counting measure.[23] The characteristic function of $L_t$ is given by the Lévy-Khintchine formula

$$\mathbb{E} e^{i\xi L_t} = e^{t\psi(\xi)}, \xi \in R,$$

where $\psi$ is the Lévy symbol.

$$\psi(\xi) = ib\xi - \frac{1}{2}\sigma^2 \xi^2 + \int_{R\setminus 0} \left(e^{i\xi y} - 1 - i\xi y \chi_{|y|\leq 1}\right) \nu(dy).$$

The process that we considered in the present work is the one with the triplet $(0, 0, \nu)$ and

$$\nu(dy) = e^{-\lambda |y|^\alpha} dy.$$

## Data Availability Statement

The data that support the findings of this study are openly available in GitHub.[36]